# An Etruscan Dodecahedron

**Amelia Carolina Sparavigna**
Department of Applied Science and Technology
Politecnico di Torino, C.so Duca degli Abruzzi 24, Torino, Italy

The paper is proposing a short discussion on the ancient knowledge of Platonic solids, in particular, by Italic people.

How old is the knowledge of Platonic solids? Were they already known to the ancients, before Plato? If we consider Wikipedia [1], the item on Platonic solids is telling that there are some objects, created by the late Neolithic people, which can be considered as evidence of knowledge of these solids. It seems therefore that it was known, may be a millennium before Plato, that there were exactly five and only five perfect bodies. These perfect bodies are the regular tetrahedron, cube, octahedron, dodecahedron and icosahedron.

In his book on regular polytopes [2], Harold Scott Macdonald Coxeter, writes "The early history of these polyhedra is lost in the shadows of antiquity. To ask who first constructed them is almost as futile to ask who first used fire. The tetrahedron, cube and octahedron occur in nature as crystals. ... The two more complicated regular solids cannot form crystals, but need the spark of life for their natural occurrence. Haeckel (Ernst Haeckel's 1904, Kunstformen der Natur.) observed them as skeletons of microscopic sea animals called radiolaria, the most perfect examples being the Circogonia icosaedra and Circorrhegma dodecahedra. Turning now to mankind, excavations on Monte Loffa, near Padua, have revealed an Etruscan dodecahedron which shows that this figure was enjoyed as a toy at least 2500 years ago."

Before Plato, Timaeus of Locri, a philosopher among the earliest Pythagoreans, invented a mystical correspondence between the four easily constructed solids (tetrahedron, icosahedron, octahedron and cube), and the four natural elements (fire, air, water and earth). "Undeterred by the occurrence of a fifth solid, he regarded the dodecahedron as a shape that envelops the whole universe." [2].

It is interesting that Donald Coxeter is reporting the existence of an Etruscan dodecahedron, that is, an object having the shape of a Platonic solid found in Italy, not of Greek origin.  In Refs.3 and 4 too, it is told that there exists an Etruscan dodecahedron made of soapstone found near Padua and believed to date from before 500 BC. Another book referring to this dodecahedron is Ref.5, is that written by György Darvas.

György Darvas discusses in [5] the Platonic solids and their use as dice. He tells that the best known of them is the cube. We use it in gambling, "because of its symmetries, it is equally likely to fall on any of its sides. … In truth, any regular body satisfies this condition of falling on any side with the same probability, not just the six-sided cube, that we in contemporary Europe are accustomed to call dice in this context.". The author continues telling that etymologically, "the noun *dice* does not even refer to a cube. This is the plural of the noun *die,* here meaning a surface with a relieved design forming one of the facets of polyhedron.

In principle, any of the five regular polyhedra can be used as a die. "There is an evidence to suggest that in Italy of old, dodecahedra were used in games, while in Etruscan cultures, they can have a religious significance (Figure 6.9a)". This is what reference 5 is telling.
In fact, this figure 6.9a of Ref.5 (Fig.1 shows a snapshot of what we can see by means of Google Books) is showing a Roman dodecahedron, not an Etruscan dodecahedron as the caption is telling. The book continues: "In Japan, for example - where the number five is considered a lucky mascot - a dodecahedron delimited by regular pentagons is still used for this purpose to this day. Sometimes it is customary to write the digits from one to twelve on its faces, sometimes the names of the twelve months.".

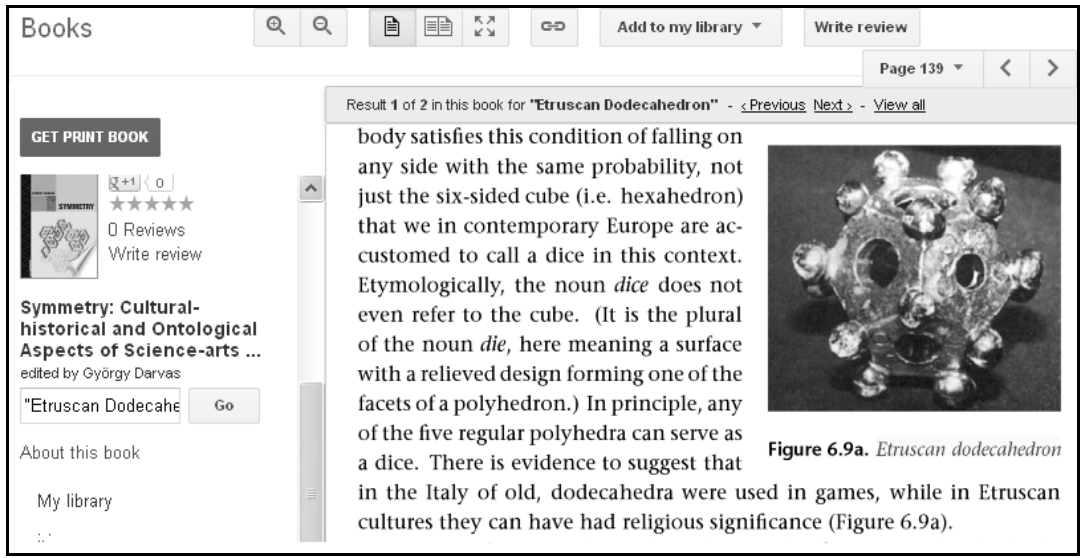

Fig.1 The image shows a snapshot of a page of Ref.4 that we can see by means of Google Book.

Fig.1 shows that Figure 6.9a of Ref.[4] can be misleading. This is a Roman dodecahedron of the second or third century AC (see Ref.6), having probably a use quite different from that of dice.
What was then the shape of the Etruscan dodecahedron? Let us report the original discussion and illustration of the researcher that found it. He was Stefano De' Stefani. In the Proceedings of the Royal Venetian Institute of Sciences, Arts and Letters ([7], 1885, see Appendix A for a small part of the Italian text), the author tells where the dodecahedron was found and reports about the existence of an icosahedron in Turin. The paper is entitled "On an almost regular dodecahedron of stone, with pentagonal faces carved with figures, discovered in the ancient stone huts of Monte Loffa".
The place of discovery belongs to Sant'Anna del Faedo village of Breonio, in the region of the western Lessini Mountains, called by the ancient historians as the region of Reti and Euganei, who were destroyed and scattered by the Gauls. De' Stefani is in agreement with several ancient writers, who considered Reti an ancient Italic people of Etruscan origin, that under the Gauls pressure had to find refuge on Alps [8].
The author continues telling that Gauls, "people of wild and fierce aspect", leaved in the same huts of Monte Loffa the manifest evidence of their presence, shown by tools, weapons and ornaments. "This village or encampment of prehistoric times shows objects of human industry that are represented by flint tools and weapons from the Neolithic period, of Etruscan bronzes type or Euganeo and Gaulish coins and other objects.
The paper has an illustration showing the dodecahedron (see Fig.2).

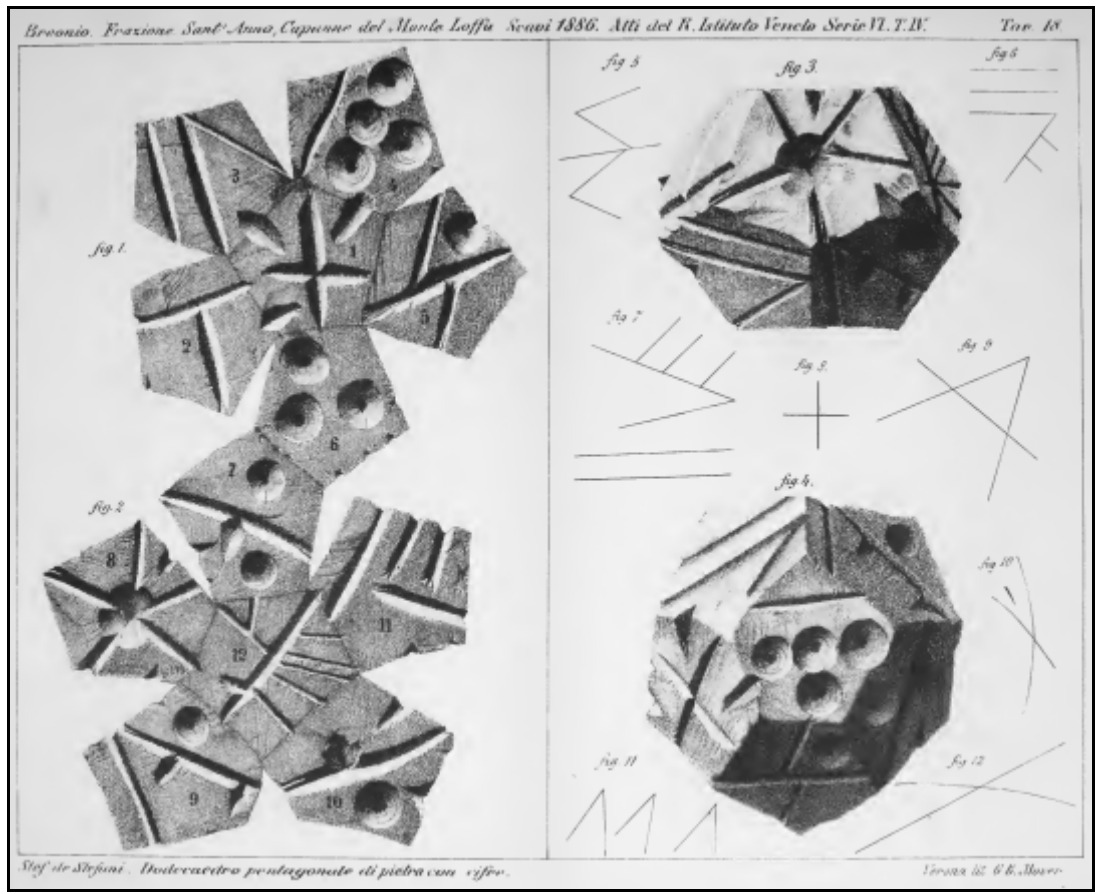

Fig.2 Etruscan dodecahedron from Monte Loffa (from ref.7).

The paper continues with a deep discussion of the nature and use of the dodecahedron in Fig.2. Several scholars were interviewed by De' Stefani, and he came to the conclusion that this dodecahedron was a die.
The paper [7] reports the opinion of Ariodante Fabretti [9], that De' Stefani received in a letter written by Carlo Cipolla [10]. Fabretti says that it is a die. The signs are conventional, perhaps a sort of numerals. In this case, this specimen is interesting because it seems to show a mixture of dots, as in our modern dice, and Etruscan numbers, adapted from the Greek numerals. On one of the face we can see "IV", may be, for "four".
There is also another interesting fact. Fabretti showed to Cipolla an icosahedron that could had some link with this dodecahedron. The icosahedron was made of blue-glazed earthenware. On each face there were impressed some Greek letters. Cipolla asked Fabretti if he knew anything about the origin of the icosahedron. He replied that it was owned by the city of Turin, before coming to the Museum of Antiquities, on occasion of an exchange. It was therefore supposed that this object was found in Piedmont.

It seems that in 1885, the existence of the icosahedron was unpublished. Unfortunately, I do not know where the Turin icosahedron is. Probably it is like that shown in Fig.3, from the second century AD, sold at auction for about $18,000 [11]. In my opinion, the Turin icosahedron could be older.

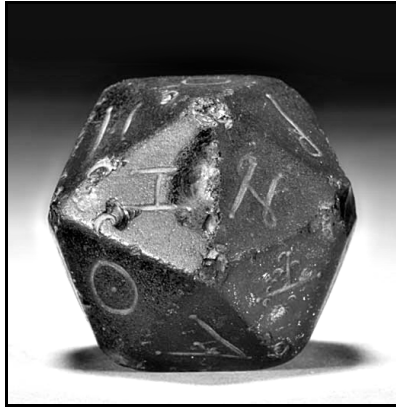 Fig3. The icosahedron die of Ref.11.

We could conclude that the ancient people in Italy, trading with Greeks, imported some numerals, and, among the first applications, used them on dice for gambling. In any case, they developed their own numeral system that evolved in the Roman numeral system.

**Appendix A**
Here the reader can find a small part, in Italian, of the paper "Intorno un dodecaedro quasi regolare di pietra a facce pentagonali scolpite con cifre, scoperto nelle antichissime capanne di pietra del Monte Loffa. by Stefano De' Stefani, Atti del Reale Istituto veneto di scienze, lettere ed arti (1885). The author tells where the dodecahedron was found and reports of the existence of an icosahedron in Turin.

"Il luogo poi del rinvenimento del nostro dado appartiene a Sant'Anna del Faedo o d'Alfaedo, frazione del Comune di Breonio, nella parte dei monti Lessini occidentali, chiamata dagli antichi storici regione dei Reti e degli Euganei, i quali sarebbero stati poscia rotti e dispersi dai Galli. Questo popolo feroce di selvaggio aspetto, dedito alle rapine ed alle stragi, avrebbe lasciate nelle stesse capanne del Monte Loffa le prove manifeste della sua presenza e della sua barbarie. La sua presenza sarebbe dimostrata da arnesi, armi ed ornamenti di tipo generalmente chiamato gallico..., che io stesso raccolsi fra i carboni o le macerie accumulate sotto le grandi lastre di pietra del luogo, che formavano il tetto delle capanne. Questo villaggio od accampamento dei tempi protostorici, mostra oggetti dell'umana industria: Ivi sono rappresentati da armi ed arnesi di selce del periodo neolitico, da bronzi di tipo etrusco od euganeo (Esle e Certosa) e da oggetti diversi e monete galliche…"

L'autore riporta poi il parere di Ariodante Fabretti sul dodecaedro, che gli è stato comunicato con lettera da Carlo Cipolla. Cipolla riferisce che Fabretti dice che trattasi di un oggetto "lusorio", una specie di dado. I segni sarebbero segni onvenzionali, una specie di segni numerali forse. Il Fabretti mostrò a Cipolla un icosaedro che aveva delle attinenze col suo dodecaedro. Questo icosaedro era di una pasta terrosa, smaltata. Dello smalto restavano molte parti; di un bellissimo colore cilestro. Sopra ciascuna faccia si leggeva impressa una chiara lettera greca. Secondo il Fabretti questo icosaedro avrebbe appunto servito per giuoco. Cipolla chiese al prof. Fabretti s'egli avesse notizie sulla provenienza del suo icosaedro. Egli rispose che apparteneva in addietro al Municipio di Torino, e che provenne al Museo di Antichità, in occasione di un cambio fatto. Non dubita che al Municipio sia stato dato da B. Gastaldi; quindi è più che probabile che detto oggetto sia stato rinvenuto nel Piemonte. L'icosaedro, dice. Fabretti a Cipolla, è tuttora inedito.